\documentclass{gtmon_a}
\pdfoutput=1
\usepackage{amscd}
\usepackage[matrix,arrow,curve]{xy}

\proceedingstitle{Proceedings of the Nishida Fest (Kinosaki 2003)}
\conferencestart{28 July 2003}
\conferenceend{8 August 2003}
\conferencename{International Conference in Homotopy Theory}
\conferencelocation{Kinosaki, Japan}

\editor{Matthew Ando}
\givenname{Matthew}
\surname{Ando}

\editor{Norihiko Minami}
\givenname{Norihiko}
\surname{Minami}

\editor{Jack Morava}
\givenname{Jack}
\surname{Morava}

\editor{W Stephen Wilson}
\givenname{W Stephen}
\surname{Wilson}

\title{On fibrations related to real spectra}

\author{Nitu Kitchloo}
\givenname{Nitu}
\surname{Kitchloo}
\address{Department of Mathematics\\
University of California, San Diego\\\newline
La Jolla CA 92093-0112\\
USA}
\email{nitu@math.ucsd.edu}
\urladdr{}

\author{W Stephen Wilson}
\givenname{W Stephen}
\surname{Wilson}
\address{Department of Mathematics\\
Johns Hopkins University\\\newline
Baltimore MD 21218\\
USA}
\email{wsw@math.jhu.edu}
\urladdr{}

\volumenumber{10}
\issuenumber{}
\publicationyear{2007}
\papernumber{14}
\startpage{237}
\endpage{244}

\doi{}
\MR{}
\Zbl{}

\arxivreference{}


\received{24 June 2004}
\revised{11 September 2005}
\accepted{}
\published{27 January 2007}
\publishedonline{27 January 2007}
\proposed{}
\seconded{}
\corresponding{}
\version{}

\makeatletter
\def\cnewtheorem#1[#2]#3{\newtheorem{#1}{#3}[section]
\expandafter\let\csname c@#1\endcsname\c@thm}

\let\xysavmatrix\xymatrix
\def\xymatrix{\disablesubscriptcorrection\xysavmatrix}
\AtBeginDocument{\let\bar\wbar\let\tilde\wtilde\let\hat\what}

\newtheorem{thm}{Theorem}[section]
\cnewtheorem{corr}[thm]{Corollary}
\cnewtheorem{lemma}[thm]{Lemma}
\cnewtheorem{conj}[thm]{Conjecture}
\cnewtheorem{prop}[thm]{Proposition}
\cnewtheorem{claim}[thm]{Claim}
\cnewtheorem{prob}[thm]{Question}
\cnewtheorem{hyp}[thm]{Hypothesis}

\theoremstyle{definition}
\cnewtheorem{defn}[thm]{Definition}
\cnewtheorem{exm}[thm]{Example}

\theoremstyle{remark}
\cnewtheorem{remark}[thm]{Remark}
\cnewtheorem{note}[thm]{Note}

\makeatother  %  move after \newtheorem block

\newcommand{\Smash} {\wedge}

\newcommand{\Wedge} {\vee}

\newcommand{\field}[1] {\mathbb #1}

\newcommand{\E} {\field E}
\newcommand{\M} {\field M\field U}

\newcommand{\B} {\field B}
\renewcommand{\D} {\field P}

\newfont{\german}{eufm10}

\newcommand{\A}{\ifmmode{\mathcal{K}}\else${\mathcal{K}}$\fi}
\newcommand{\U}{\ifmmode{\mathcal{U}}\else${\mathcal{U}}$\fi}

\DeclareMathOperator*{\dirlim}{\varinjlim}
\DeclareMathOperator*{\colim}{\colim\,}

\DeclareMathOperator{\map}{Map}
\newcommand{\ra}{\rightarrow}
\newcommand{\lra}{\longrightarrow}

\newcommand{\llra}[1]{\stackrel{#1}{\lra}}

\begin{document}

\begin{htmlabstract}
We consider real spectra, collections of <b>Z</b>/(2)&ndash;spaces indexed
over <b>Z</b>&oplus;<b>Z</b>&alpha; with compatibility conditions.
We produce fibrations connecting the homotopy fixed points and
the spaces in these spectra.  We also evaluate the map which is the
analogue of the forgetful functor from complex to reals composed with
complexification.  Our first fibration is used to connect the real
2<sup>n+2</sup>(2<sup>n</sup>-1)&ndash;periodic Johnson&ndash;Wilson
spectrum ER(n) to the usual 2(2<sup>n</sup>-1)&ndash;periodic
Johnson&ndash;Wilson spectrum, E(n).  Our main result is the fibration
&Sigma;<sup>&lambda;(n)</sup>ER(n)&rarr;ER(n)&rarr;E(n), where
&lambda;(n)=2<sup>2n+1</sup>-2<sup>n+2</sup>+1.
\end{htmlabstract}

\begin{abstract}
We consider real spectra, collections of $\mathbb{Z}/(2)$--spaces indexed
over $\mathbb{Z}\oplus\mathbb{Z}\alpha$ with compatibility conditions.
We produce fibrations connecting the homotopy fixed points and the spaces
in these spectra.  We also evaluate the map which is the analogue of the
forgetful functor from complex to reals composed with complexification.
Our first fibration is used to connect the real $2^{n+2}(2^n-1)$--periodic
Johnson--Wilson spectrum $ER(n)$ to the usual $2(2^n-1)$--periodic
Johnson--Wilson spectrum, $E(n)$.  Our main result is the fibration
$\Sigma^{\lambda(n)} ER(n) \to ER(n) \to E(n)$, where $\lambda(n) =
 2^{2n+1}-2^{n+2}+1$.
\end{abstract}

\maketitle

\section{Introduction}

In 1968, Landweber \cite{Land:Conj} introduced the idea of a real complex
cobordism by taking the homotopy fixed points of complex cobordism
under complex conjugation.  A few years later this theory was studied
again by Araki and Murayama \cite{Araki:forget,Araki:orient,Araki:real}.
Recently there has been a flurry of activity around this theory by Hu
and Kriz \cite{HK5,HK6,HK,HK3,HK2,HK4}.

In \cite{HK}, Hu and Kriz  produce real versions, $ER(n)$,
of the Johnson--Wilson spectra $E(n)$ (see Johnson and Wilson
\cite{JW2}) and compute their homotopy.  The homotopy of $E(n)$
is $\Z_{(2)}[v_1,v_2,\ldots,v_n^{\pm 1}]$.  $E(n)$ is periodic
of period $|v_n| = 2(2^n-1)$, $ER(n)$ is periodic of period
$\bigl|\smash{v_n^{2^{n+1}}}\bigr| = 2^{n+2}(2^n-1)$, and the construction
gives maps of spectra: $ER(n) \ra E(n)$.

In the case $n=1$ this is just the map $KO_{(2)} \ra KU_{(2)}$ and
Wood identified the fibre as $\Sigma KO_{(2)}$.  

The main purpose of this paper is to identify the fibre of
$ER(n) \ra E(n)$, producing the fibration:
\begin{equation}
\label{firstfib}
\Sigma^{\lambda(n)} ER(n) \ra ER(n) \ra E(n).
\end{equation}
where $\lambda(n)= 2^{2n+1}-2^{n+2}+1$.
These fibrations should make these theories much more accessible.

Let $\E$ be a real spectrum as defined in \cite{HK}. In particular, $\E$ 
is given by a
collection of pointed $\Z/(2)$--spaces $\E_V$ indexed by the 
representation ring
$RO(\Z/(2))$ of the group $\Z/(2)$.
Recall that $RO(\Z/(2)) = \Z \oplus \Z\alpha$, where $\alpha$ is the sign 
representation. Moreover, we require that the spaces $\E_V$ be compatible 
in the following sense.

Given a representation $U$, and a pointed $\Z/(2)$--space $X$, let 
$\Omega^UX$ denote the
space $\map_*(S^U,X)$, where $S^U$ is the one-point compactification of 
$U$. The space
$\Omega^UX$ has an induced diagonal action of the group $\Z/(2)$. For the 
spectrum $\E$, we require the existence of a family of equivariant 
homeomorphisms $\alpha_{U,V} \co \Omega^U\E_{U\oplus V} \rightarrow \E_V$, 
that satisfy obvious compatibility. 

A multiplicative real spectrum $\E$ is one that admits a multiplication 
preserving the real structure (see \cite{HK}). 

\begin{exm} The real complex bordism spectrum $\M$ is defined as follows. 
Let $MU(n)$
denote the Thom space of the universal bundle over $BU(n)$. Complex 
conjugation induces an
action of $\Z/(2)$ on $MU(n)$. Define $\M_V$ as the space $
\dirlim_n\Omega^{n(1+\alpha)-V}MU(n)$ for $V \in RO(\Z/(2))$. Notice that 
$n(1+\alpha)-V$
is a well defined representation of $\Z/(2)$ for sufficiently large values 
of $n$. It is left to the reader to verify that $\M$ has the properties of 
a multiplicative real spectrum. 
\end{exm}

\begin{exm} The Brown--Peterson spectrum has a real analogue $\B\D$. The 
real Johnson--Wilson spectra $\E (n)$ may also be defined along similar 
lines \cite{HK}. These spectra are 
in fact multiplicative real spectra. $\E(1)$ is 2--localized real K--theory 
of Atiyah \cite{HK}.
\end{exm}

We will use the notation $\E\R_V$ to denote the homotopy fixed points of 
the
$\Z/(2)$--action on $\E_V$. Notice that for a fixed $V \in RO(\Z/(2))$, 
the collection of spaces $\{\E\R_{n+V}, n \in \Z\}$ form a spectrum in the 
usual sense. We shall abuse notation and refer to the spectra 
$\{\E\R_{n+V}, n \in \Z\}$ and $\{\E_{n+V}, n \in \Z\}$ as the spectra 
$\E\R_V$ and $\E_V$ respectively. The purpose of this paper is to relate 
$\E\R_V$ to $\E_V$ via a fibration. Of particular interest to us will be 
the case when the spectrum $\E$ is $\E(n)$. 
We need the following result which we assume is well known to the experts:

\begin{prop} \label{fibrations}
There are fibrations of spectra:
\[ \E\R_{V-\alpha} \llra{a} \E\R_V \llra{\iota} \E_V, \quad \quad \E_V 
\llra{1+\sigma} \E\R_V \llra{a} \E\R_{V+\alpha}, \]
where the map $a$ is induced by the map $a \co S^0 \lra S^{\alpha}$ given by 
the inclusion of the poles. The map $\iota$ is the standard inclusion, and 
the map $(1+\sigma)$ is a lift of the Norm map on $\E_V$. Moreover, if 
$\E$ is a multiplicative real spectrum, then $\E\R_V$ is a 
$\E\R_0$--module spectrum for all $V$, and the above fibrations are 
fibrations of $\E\R_0$--module spectra. 
\end{prop}

\begin{remark}
On the level of individual spaces we have fibrations
$$
\E\R_{m+(n-1)\alpha} \ra
\E\R_{m+n\alpha} \ra
\E_{m+n\alpha}. 
$$
This is a great help to computations as we hope to demonstrate in a future
paper.
\end{remark}

Observe that the spaces $\E_{V-1} = \Omega \E_V$ 
and $\E_{V-\alpha}=\Omega^{\alpha} \E_V$ are homeomorphic. 
(Actually, $\E_{m+n\alpha}$
and
$\E_{m'+n'\alpha}$
are the same when $m+n = m' + n'$.)
In the statement of the next theorem, we will use this 
homeomorphism to identify the two spaces. Note, however that the 
action of $\Z/(2)$ on the two spaces is different. If we let $\sigma$ 
denote the action of the generator of $\Z/(2)$ on $\E_{V-1}$, 
and $\tilde{\sigma}$ the action on $\E_{V-\alpha}$, then the 
two actions are related via $\tilde{\sigma} = \sigma \alpha_* = -\sigma$. 

Now consider the boundary map. This map $\partial$ is defined as the 
map $ \E_{V-1} \ra \E\R_{V-\alpha}$ given by looping back the first 
fibration above composed with the map $\E\R_{V-\alpha} \ra \E_{V-\alpha}$ 
given by the inclusion of the fixed points. Therefore
\[ \partial \co \E_{V-1} \lra \E_{V-\alpha}. \]
We have the following proposition.

\begin{prop}\label{boundary}
Let $\E_{V-1}$ be identified with the space $\E_{V-\alpha}$ as explained 
above. 
Then the map $\partial$ is given by $\partial = Id - \sigma = Id + 
\tilde{\sigma}$. 
\end{prop}

The standard example of this result is the composition $KU \ra KO \ra KU$
and this is just a generalization of it.  The boundary is the composition 
of
two maps.  The first can be thought of as forgetting the complex structure
and looking only at the underlying real structure.  The next map can be
thought of as complexification.  This boundary map comes in useful in 
calculations
we hope will appear in a future paper.

Our primary interest is the case when $\E = \E (n)$. 
The following theorem uses the computation of the homotopy of $\E\R (n)$ 
given in \cite{HK}.

\begin{thm}\label{ER(n)fib}
There exist nontrivial elements $x(n) \in \pi_{\lambda(n)}(\E\R(n)_0)$, 
where $\lambda(n)$ is the integer defined by $\lambda(n) = 
2^{2n+1}-2^{n+2}+1$, such that one has a fibration of $\E\R(n)_0$--module 
spectra:
\[ \Sigma^{\lambda(n)}\E\R(n)_V \llra{ x(n)} \E\R(n)_V \llra{\iota} 
\E(n)_V. \]
\end{thm}

\begin{remark}
An interesting special case of the above theorem is when $n=1$, and $V=0$. 
Note that $\E(1) = KU_{(2)}$, and hence $\E\R(1) = KO_{(2)}$. Moreover, 
the element $x(1)$ is none other than $\eta$.
Hence one reproduces a well-known result
\[ \Sigma KO_{(2)} \llra{ \eta} KO_{(2)} \lra KU_{(2)}. \]
More generally, fixing $V=0$, we get the fibration \eqref{firstfib}.
\end{remark}

Our dependence on the published work of Hu and Kriz is obvious.  
In addition, they recently informed us
they can prove a generalization of our result.  
By working with the universal
example, ie inverting $v_n$ in $MU$, they can show that any theory with
$v_n$ inverted has the same fibration we have for $ER(n)$.  
The proof
is the same.
Note that this works for their version of real Morava $K$--theory in
\cite{HK} making it a much more interesting theory but unfortunately
still not a ring theory.

\section{The fibrations}

In this section we will show the existence of the two fibration given in 
the introduction. 

Let $S^{\alpha}$ denote the one-point compactification of the one 
dimensional nontrivial
representation of $\Z/(2)$. Notice that one has a $\Z/(2)$--equivariant 
cofibration:
\begin{equation} \label{cof} \Z/(2)^+ \lra S^0 \llra{a} S^{\alpha} 
\end{equation}
where the map $\Z/(2)^+ \ra S^0$ is given by the pinch map. Let $\E$ be a 
real spectrum, and for the purposes of this section, let $\E_V$ denote the 
spectrum given by the collection of spaces $\{\E_{n+ V}, n\in \Z \}$. 
Smashing the cofibration \eqref{cof} yields a cofibration of equivariant 
spectra
\begin{equation} \label{Ecof} \E_V \Smash \Z/(2)^+ \lra \E_V \llra{a} 
\E_{V+\alpha}. \end{equation}
Notice that $\Z/(2)^+$ may be identified with $S^0 \Wedge S^0$, with the 
$\Z/(2)$ action
given by the twist map. Under this identification, the pinch map $\Z/(2)^+ 
\ra S^0$
corresponds to the fold map $S^0 \Wedge S^0 \ra S^0$. Hence, in the 
category of spectra,
$\E_V \Smash \Z/(2)^+$ may be identified with $\E_V \Wedge \E_V = \E_V 
\times \E_V$ with
the $\Z/(2)$ action given by $\tilde{\sigma}(x,y) = 
(\sigma(y),\sigma(x))$, where $\sigma$
denotes the generator of $\Z/(2)$.
Furthermore, the pinch map $\E_V \Smash \Z/(2)^+ \ra \E_V$ corresponds to 
the sum map $\E_V \times \E_V \llra{+} \E_V$. 

Consider the twisted diagonal map 
\[ \Delta \co \E_V \lra \E_V \times \E_V, \quad \quad \Delta(x) = 
(x,\sigma(x)). \]
Notice that $\tilde{\sigma}\Delta(x) = \Delta(x)$. From this it follows 
easily that
$\Delta$ lifts to an equivalence $\E_V \ra (\E_V \times \E_V)^{h\Z/(2)}$. 
Putting these results together, we get the following proofs.

\begin{proof}[Proof of the second fibration in \fullref{fibrations}.]
Taking homotopy fixed points of \ref{Ecof} yields another fibration. If we 
identify $(\E_V
\Smash \Z/(2)^+)^{h\Z/(2)}$ with $\E_V$, then the map $(\E_V \Smash 
\Z/(2)^+)^{h\Z/(2)}
\ra (\E_V)^{h\Z/(2)}$ is a lift of $(1+\sigma)$. 
\end{proof}

\begin{proof}[Proof of the first fibration in \fullref{fibrations}.]
For the second fibration, one considers the Spanier--Whitehead dual of 
\ref{cof}:
\[ S^{-\alpha} \llra{a} S^0 \lra \Z/(2)^+ \]
where the map $S^0 \ra \Z/(2)^+ = S^0 \Wedge S^0$ corresponds to the 
diagonal. Smashing with $\E_V$ yields an equivariant fibration
\[ \E_{V-\alpha} \llra{a} \E_V \lra \E_V \times \E_V. \] 
Taking homotopy fixed points of this fibration and making the 
identifications described earlier, we get the remaining fibration:
\[ \E\R_{V-\alpha} \llra{a} \E\R_V \llra{\iota} \E_V. \]
To complete the proof one simply observes that all the above constructions 
respect the $\E\R_0$--module structure if $\E$ is a multiplicative real 
spectrum.
\end{proof}

\section{The boundary map}
In this section, we analyse the boundary map for the above fibrations. 
This map $\partial$ is defined as the composite of the map $ \E_{V-1} \ra 
\E\R_{V-\alpha}$ given by looping back the fibration constructed in the 
previous section, and the map $\E\R_{V-\alpha} \ra \E_{V-\alpha}$ given by 
the inclusion of the fixed points. Therefore
\[ \partial \co \E_{V-1} \lra \E_{V-\alpha}. \]
The map $\partial$ may be explicitly constructed as follows. Consider the 
composite equivariant map given by the fold map followed by the pinch map:
\begin{equation}\label{fib2} \Z/(2)^+ \Smash S^{\alpha} \llra{f} 
S^{\alpha} \llra{p} \Z/(2)^+ \Smash S^1. 
\end{equation}
Notice that the Spanier-Whitehead dual of the pinch map $ p \co S^{\alpha} 
\ra \Z/(2)^+ \Smash S^1$ is the difference map $ (-) \co Z/(2)^+ \Smash 
S^{-1} \ra S^{-\alpha}$. Taking the Spanier-Whitehead dual of the 
composite \eqref{fib2} yields 
\[ \Z/(2)^+ \Smash S^{-1} \llra{(-)} S^{-\alpha} \llra{\Delta} 
\Z/(2)^+\Smash S^{-\alpha}. \]
On smashing the above with $\E_V$, we obtain the composite map
\[ \Delta(-) \co \E_{V-1} \times \E_{V-1} \llra{(-)} \E_{V-\alpha} 
\llra{\Delta} \E_{V-\alpha}
\times \E_{V-\alpha}. \]
From the previous section, we can see that there is a commutative diagram
$$
\xymatrix{
 \E_{V-1} \ar[r]^{\partial} \ar[d]^{(1,\sigma)} & \E_{V-\alpha} 
\ar[d]^{(1,\tilde{\sigma})} \\
\E_{V-1} \times \E_{V-1} \ar[r]^{\Delta(-)} & \E_{V-\alpha} \times 
\E_{V-\alpha} }
$$
where $\sigma$ denotes the $\Z/(2)$--action on $\E_{V-1}$, and 
$\tilde{\sigma}$ denotes the
$\Z/(2)$--action on $\E_{V-\alpha}$. Recall that the spaces $\E_{V-1} = 
\Omega \E_V$ and $\E_{V-\alpha} = \Omega^{\alpha}\E_V$ are homeomorphic 
and the above two actions are related via $\tilde{\sigma} = 
\sigma\alpha_*$. Since $\alpha_*$ is homotopic to the inversion, we have 
$\tilde{\sigma} = -\sigma$. From a diagram chase we get $\partial (x) = x 
- \sigma(x) = x + \tilde{\sigma}(x)$. 

\section[The case of $E(n)$]{The case of $\E(n)$}

We recall the computation (via the Borel spectral sequence) of the 
homotopy of $\B\D\R$ given in \cite{HK}, and described in the form we need 
in \cite{HK5}. We will reproduce the Borel spectral sequence with $\B\D$ 
replaced by $\E(n)$. The $E_2$--term of the Borel spectral sequence for 
$\E(n)$ is given by 
\[ E_2 = \Z_{(2)}[v_k, v_n^{\pm 1}, a, \sigma^{\pm 2}]/(2a), \qua n > k 
\geq 0, \qua v_0
= 2. \]
The bidegrees of the generators are given by 
\[ |a| = -\alpha, \qua |v_k| = (2^k-1)(1+\alpha), \qua |\sigma^2| = 
2(\alpha-1). \]
The differentials are given by comparing with the Borel spectral sequence 
converging to the homotopy of $\B\D\R$. In particular, the elements $v_k$ 
and $a$ are permanent cycles, and the nontrivial differentials are 
\[ d_{2^{k+1}-1}(\sigma^{-2^k}) = v_ka^{2^{k+1}-1}, \quad 0 < k \leq n. \]
Using the methods of \cite{HK}, \cite{HK5}, we notice that the 
$E_{\infty}$--term for the homotopy of $\E\R(n)$ is given by the following 
ring:
\[ \Z_{(2)}[v_k\sigma^{l2^{k+1}}, a, v_n^{\pm 1}, \sigma^{\pm 2^{n+1}}]/I, 
\qua n > k\geq 0,
\qua l \in \Z \]
where $I$ is the ideal generated by the relations:
\begin{align*} v_0 & =2, \\
a^{2^{k+1}-1}v_k\sigma^{l2^{k+1}}& =0, \\
v_m\sigma^{l2^{m+1}}. v_k\sigma^{s2^{m-k}2^{k+1}} & = 
v_k.v_m\sigma^{(l+s)2^{m+1}}\qua
m\geq k.
\end{align*}
The bidegrees of the generators are given by
\[ |a| = -\alpha, \quad \bigl|v_k\sigma^{l2^{k+1}}\bigr| = (2^k-1)(1+\alpha) + 
l2^{k+1}(\alpha-1). \]
Comparing with the homotopy of $\B\D\R$, we notice that there are no 
extension problems, and so the above is in fact isomorphic to the 
homotopy of $\E\R(n)$.

Now consider the element 
\[ y(n) = v_n^{2^n-1}\sigma^{-2^{n+1}(2^{n-1}-1)}, \quad y(n) \in 
\pi_{\lambda(n)}(\E\R(n)_{-\alpha}) \]
where $\lambda(n)$ is the integer $\lambda(n) = 2^{2n+1}-2^{n+2}+1$. The 
element $y(n)$ is clearly invertible in the above ring. Hence we get 
the following claim.

\begin{claim} \label{twist}
Multiplication by the element $y(n)$ yields an equivalence of 
$\E\R_0$--module spectra:
\[ \Sigma^{\lambda(n)}\E\R(n)_V \llra{ y(n)} \E\R(n)_{V-\alpha}. \]
\end{claim}
We define the element $x(n)$ to be the element 
\[ x(n) = a.y(n), \quad x(n) \in \pi_{\lambda(n)}(\E\R(n)_0). \]

This claim, along with the first fibration given in \fullref{fibrations} yields 
the proof of \fullref{ER(n)fib}. 

\begin{remark}
The spectrum $\E\R(n)_0$ is periodic with period $2^{n+2}(2^n-1)$ 
generated by the homotopy
element $v_n^{2^{n+1}}\sigma^{-2^{n+1}(2^n-1)}$. 
\end{remark} 

\bibliographystyle{gtart}
\bibliography{link}

\end{document}